\documentclass[11pt,a4paper]{amsart}

\usepackage{amsmath,amssymb,latexsym,amsthm,enumerate}
\usepackage{verbatim}
\usepackage{graphicx} 
\usepackage{hyperref} 

\def\RR{\mathbb{R}}

\def\1{\mathbf{1}}
\def\0{\mathbf{0}}

\def\I{\mathbf{I}}

\newcommand{\te}{\widehat{e}}

\newtheorem{Thm}{Theorem}

\newtheorem{Lemma}{Lemma}

\newtheorem{As}{Assumption}

\def\BEN{\begin{enumerate}}  \def\BI{\begin{itemize}}
\def\EEN{\end{enumerate}}   \def\EI{\end{itemize}}
    \def\nn{\nonumber}

\def\mbb{\mathbb}  
\def\mc{\mathcal}  \def\ovl{\overline}

\def\te#1{\mathrm{e}^{#1}}   

\def\WH{\widehat}

\def\g{\gamma}     
\def\ve{\varepsilon}
\def\e{\epsilon}
 
  \def\nn{\nonumber}   
     
 \def\q{\qquad} 
   
  \def\td{\text{\rm d}}
 
\numberwithin{equation}{section}

\newcommand{\exit}{{\mbox{\, \vspace{3mm}}}
\hfill\mbox{$\square$}}

\begin{document}

\title[Joint limit law of undershoots and overshoots of reflected processes]{
Buffer-overflows: joint limit laws of undershoots
and overshoots of reflected processes}

\author{Aleksandar Mijatovi\'{c}}
\address{Department of Mathematics, Imperial College London}
\email{a.mijatovic@imperial.ac.uk}
\author{Martijn Pistorius}
\address{Department of Mathematics, Imperial College London}
\email{m.pistorius@imperial.ac.uk}


\thanks{MP's research supported by NWO-STAR}


\maketitle

\begin{abstract}
Let $\tau(x)$ be the epoch of first entry into 
the interval $(x,\infty)$, $x>0$, of the reflected process~$Y$ 
of a L\'evy process $X$, and define the overshoot $Z(x) = Y(\tau(x))-x$ and 
undershoot $z(x) = x - Y(\tau(x)-)$ of 
$Y$ at the first-passage time over the level $x$. In this paper we 
establish, separately under the Cram\'{e}r and positive drift assumptions, 
the existence of the weak limit of
$(z(x), Z(x))$ as $x$ tends to infinity 
and provide explicit formulae for their joint CDFs
in terms of the L\'{e}vy measure of $X$ and the 
renewal measure of the dual of $X$. We apply our results 
to analyse the behaviour of the classical M/G/1 queueing system at
the buffer-overflow, both in a stable and unstable case.
\end{abstract}

\section{Introduction}
Consider a classical single server M/G/1 queueing system, 
consisting of a stream of jobs of sizes given by positive IID 
random variables 
$U_1, U_2, \ldots$, 
arriving according to a standard Poisson process $N$ with rate
$\lambda$, and a server that processes these jobs at unit speed. 
At a given time $t$ the workload $Y$ in the system is given by 
\begin{eqnarray}\label{eq:refl}
&& Y(t) = X(t) - X_*(t), \qquad X_*(t) = \inf_{s\leq t} \{X(s),0\},\\
&& X(t) = I(t) - O(t),\qquad I(t) = u_0 + \sum_{n=1}^{N_t} U_n, 
\q O(t) = t, \label{eq:cpp}
\end{eqnarray}
where $u_0\ge 0$ is the workload in the system at time $0$, 
$I(t)$ denotes the cumulative
workload of all jobs that have arrived by time $t$ and 
$O(t)$ the cumulative capacity at time~$t$ (i.e. the amount 
of service that could have been provided if the server has never been 
idle up to time $t$). We refer to~\cite{AsmussenAPQ,Prabhu} 
for background on queueing theory. In generalisations of the classical 
M/G/1 model it has been proposed to replace the
compound Poisson process $X$ in Eqn.~\eqref{eq:cpp} by a general L\'{e}vy
process leading to the so-called L\'{e}vy-driven queues.
In case the system in Eqns.~\eqref{eq:refl}--\eqref{eq:cpp} has a finite
buffer of size $x>0$ for the storage of the workload, two quantities
of interest are the under- and overshoots of the workload $Y$ at the first time
$\tau(x)$ of \textit{buffer-overflow} (i.e. $z(x) = x - Y(\tau(x)-)$
and $Z(x) = Y(\tau(x))-x$ resp.),
representing the level of the workload just
before the buffer-overflow 
and the part of the job lost at $\tau(x)$
(see e.g.~\cite{Chydzinski, Kampa}). 

Our main result
(Theorem~\ref{thm} below)
states that 
if a 
L\'evy process $X$ satisfies the Cram\'er assumption
and a non-lattice condition, 
the joint weak limit 
$(z(\infty),Z(\infty))$
of the 
under- and overshoot of $Y$ at $\tau(x)$ (as $x\uparrow\infty$) 
exists and is explicitly given by the following formula:
\begin{equation}
\label{eq:Cramer_Joint_law}
P\left[(z(\infty),Z(\infty))\in\td u\otimes\td v\right] 
 = \frac{\gamma}{\phi(0)} \WH V_\gamma(u)\,
\nu(\td v+u)\,\td u,
\qquad\text{for $u,v>0$,}
\end{equation}
where
$\WH V_\gamma(u) \doteq \int_{[0,u]}\te{\gamma (u-y)} \WH V(\td y)$,
$\WH V(\td y)$ 
is the renewal measure of 
the dual of $X$,
$\phi$
the Laplace exponent of the ascending ladder-height process, 
$\nu$
the L\'evy measure of 
$X$
and 
$\gamma$ the Cram\'er coefficient (see
Section~\ref{sec:main} for precise definitions). 
Note that, unlike~\eqref{eq:Cramer_Joint_law}, 
analogous limit results for the L\'evy process
$X$ 
at the moment of its first passage over the level
$x$
require conditioning on the event that the process 
reaches the level 
$x$
in finites time
(see e.g. Lemma~\ref{lem:Xconv}(iii) below
and~\cite[Thm.~7.1]{Griffin_Maller}, \cite[Thm.~4.2]{KluppelbergKM}).

In the case 
$E[X(1)]$
is positive and finite
and
$X(1)$
is non-lattice,
we identify (in Theorem~\ref{thm:P} below)
the limiting joint 
law 
of the under- and overshoot of $Y$
as:
\begin{equation}
\label{eq:pos_drift_joint_law}
P\left[(z(\infty),Z(\infty))\in\td u\otimes\td v\right] 
 =\frac{\WH V(u)}{\phi'(0)}\, \nu(\td v+u)\,\td u,  
\qquad\text{for $u,v>0$.}
\end{equation}

In the classical queueing model
given by~\eqref{eq:refl}--\eqref{eq:cpp}, 
the random variable
$X(1)$
is non-lattice if 
the distribution $F$ of the job-sizes 
$U_i$
is.
Our results yield the joint limit law explicitly
in terms of the distribution
$F$ and the 
arrival rate $\lambda$: 
\begin{equation}
\label{eq:Joint_LAw_Queue}
P[(z(\infty), Z(\infty))\in \td u\otimes\td v] = \frac{\lambda}{\lambda m - 1}(1 - \te{r^* u}) F(\td v + u)\td u,
\q u, v>0,
\end{equation}
where $m\doteq\int_{(0,\infty)}x F(\td x)$ is the mean of $F$ and 
$r^*$ denotes the largest 
(resp. smallest) root $s$ in
$\RR$
of the characteristic equation
$$
\lambda \int_{(0,\infty)} \te{sy} F(\td y)  - \lambda - s = 0,
$$
in the Cram\'{e}r (resp. positive drift) case,\footnote{Note that in the Cram\'{e}r case it holds  
$r^*>0$,
$\lambda m<1$
and in the (unstable) positive drift case we have
$r^*<0$,
$\lambda m>1$.}
cf. Remarks~(i) and~(ii) after Theorem~\ref{thm}
in the next section.


While the formulae in~\eqref{eq:Cramer_Joint_law},~\eqref{eq:pos_drift_joint_law}  
and~\eqref{eq:Joint_LAw_Queue}  hold for any starting point
$Y_0=u_0\ge 0$, it is not hard to see that it suffices to establish those
relations just for the value $u_0=0$.
In the Cram\`er case, the probability that the first time of buffer-overflow
over the level $x$ occurs {\em before} the end of the {\em busy period} (i.e.
before the first time that $Y$ reaches zero) tends to zero as $x$ tends to
infinity. Hence we may assume, by the strong Markov property of
$Y$ applied at the end time of the busy period,
that we have $Y_0=0$. In the positive
drift case
(i.e. when the queueing system is \textit{unstable}),
the probability of the complementary event that the buffer-overflow of size $x$
occurs {\em after} the first visit of $Y$ to the origin tends to zero 
(see Section~\ref{subsec:Proof_of_thm_1}). Hence, the proof of Theorem~\ref{thm:P} yields that 
the joint limit 
law in~\eqref{eq:pos_drift_joint_law}  is equal to the asymptotic
distribution, as $x\uparrow \infty$,
of the under- and overshoot of $X$ (with $X_0=u_0$) at the epoch of its first entrance into the set
$(x,\infty)$.  The latter limit law clearly does not depend on the starting level $u_0$,
due to the spatial homogeneity of the process $X$.


The arguments outlined in the previous paragraph also imply that the limit
distribution in~\eqref{eq:Cramer_Joint_law} 
(and hence~\eqref{eq:Joint_LAw_Queue}) remains valid if $Y$ is in its {\em steady state}, 
i.e. the workload process $Y$ was started according to its stationary distribution, which 
exists since, under the Cram\'{e}r assumption, $Y$ is an ergodic strong Markov process and
the corresponding queueing system is \textit{stable}.
Furthermore, it is worth noting that, in the Cram\'{e}r case, the right-hand
side of~\eqref{eq:Cramer_Joint_law} (and hence that of~\eqref{eq:Joint_LAw_Queue}) 
is in fact also equal to the asymptotic distribution of the 
under- and overshoot {\em conditional} on the event that the buffer-overflow takes
place in the busy period (this result follows directly by combining the proof
of Theorem~\ref{thm} below with the two-sided Cram\'{e}r estimate for $X$, see e.g.~\cite[Prop.~7]{MP2}).

Various aspects of the law of the reflected process have been studied recently in a
number of papers. The exact asymptotic decay of the distribution of the maximum
of an excursion, under the It\^{o}-excursion measure, was identified
in~\cite{DoneyMaller} under the Cram\'{e}r condition. Also in the Cram\'{e}r
case, the joint asymptotic
distributions of the overshoot, the maximum and the
current value of the reflected process were obtained in~\cite{MP2}.
In special cases a number of papers are devoted to the characterisation
of the law of the reflected process at the moment of buffer-overflow. 
For example, in the case of spectrally negative L\'{e}vy processes, 
the joint Laplace-transform of the pair $(\tau(x), Y_{\tau(x)})$ was obtained in 
\cite{AKP}. A sex-tuple law extension of this result, centred around the epoch of the
first-passage of the reflected process, was given in~\cite{MP}. 


The remainder of the paper is organised as follows:
the main results are stated in Section~\ref{sec:main} 
and their proofs are given in Section~\ref{sec:proofs}.
Appendix~\ref{app:proof_lemma} contains the proof of
Lemma~\ref{lem:Xconv}, which plays an important role in the 
proofs of Theorems~\ref{thm:P} and~\ref{thm}
and is stated in Section~\ref{subsec:Under_Over}.

\section{Joint limiting distributions}\label{sec:main}
Let $X$ be a L\'{e}vy process, that is, a stochastic process with
independent and stationary increments and c\`adl\`ag paths, with
$X(0)=0$, and let
 $Y=\{Y(t), t\ge 0\}$
be the reflected process of $X$ at its infimum, i.e.
\begin{equation}
Y(t) \doteq X(t) - \inf_{0\leq s\leq t} X(s), \qquad t\ge 0.
\end{equation}
The process $Y$ crosses any positive level $x$ in finite time 
almost surely, that is, the moment of first-passage 
$$\tau(x)\doteq\inf\{t\ge 0: Y(t)\in(x,\infty)\}$$
is finite with probability 1, for any $x>0$. Denote by $\Psi_x$ 
the joint (complementary) distribution function  
of the pair $(z(x), Z(x))$ of
under- and overshoot of $Y$, 
\begin{eqnarray*}
&& \Psi_{x}(u,v) \doteq P(z(x) > u, Z(x) > v), \qquad u,v\geq0,\quad x>0,
\end{eqnarray*}
where we defined 
$z(x)\doteq x - Y(\tau(x)-)$ and  $Z(x) \doteq Y(\tau(x)) - x$.

Recall that the renewal function $V:\mbb R_+\to\mbb R_+$  of $X$ 
is the unique non-decreasing 
right-continuous function with the Laplace transform given by
\begin{eqnarray}\label{eq:LTV}
\int_0^\infty \te{-\theta y} V(y) \td y & = & (\theta \phi(\theta))^{-1},\qquad \text{where}\\
\phi(s) & \doteq & - \log E[\te{-s H(1)}\I_{\{1<L(\infty)\}}],\q\text{ for $s\ge 0$},
\label{eq:Def_phi}
\end{eqnarray}
$L$ denotes a local time of $X$ at its running supremum $X^*$, 
$X^*(t) \doteq \sup_{0\leq s\leq t}X(s)$, with $L(\infty)=\lim_{t\uparrow\infty}L(t)$,  and $H$  
is the ascending ladder-height process of $X$.
The corresponding  measure $V(\td y)$ 
is the potential measure of $H$, i.e. 
$V(\td y)=\int_0^\infty P(H(t)\in \td y)\td t$. Similarly 
$\WH L$, $\WH H$, $\WH \phi$
and $\WH V$ denote the local time, the ladder process, 
its characteristic exponent and the renewal function of the dual 
process
$\WH X\doteq-X$
respectively. We assume throughout the paper that 
$L$
and $\WH L$
are normalised 
in such a way that 
$-\log E[\te{\mathrm{i}\theta X(1)}]=\phi(-\mathrm{i}\theta)\WH \phi(\mathrm{i}\theta)$,
$\theta\in\RR$,
holds,
and denote by 
$\nu$ the L\'{e}vy measure of $X$ 
and by 
$\ovl\nu(a)\doteq \nu((a,\infty))$, $a>0$,
its tail function.
For the background on ladder processes and fluctuation theory 
we refer to Bertoin~\cite[Ch.~VI]{Bertoin}. 

The first limit result concerns the positive drift case:

\begin{Thm}\label{thm:P}
Let the law of 
$X(1)$
be non-lattice
and suppose $E[|X(1)|]<\infty$, $E[X(1)]\in(0,\infty)$. Then 
$\phi'(0)\in(0,\infty)$
and the limit $\Psi_\infty(u,v)\doteq\lim_{x\uparrow\infty}\Psi_x(u,v)$ exists
and is given as follows:
\begin{equation*}
\Psi_\infty(u,v) =
\frac{1}{\phi'(0)}\int_u^\infty \ovl\nu(v+z)\WH V(z)\td z,\qquad u,v\geq0.
\end{equation*}
\end{Thm}
In the negative drift case we will restrict ourselves
to  the classical Cram\'{e}r setting:\footnote{If a non-trivial L\'evy process
$X$
satisfies As.~\ref{C}, then
$E[X(1)]\in(-\infty,0)$
since
$u\mapsto E[\te{uX(1)}]$
is strictly convex on
$[0,\gamma]$.} 
\begin{As}\label{C}\rm
Suppose that the Cram\'{e}r-assumption holds, i.e. there exists a $\gamma\in(0,\infty)$ such that $E[\te{\gamma X(1)}]=1$, 
$X(1)$
is non-lattice with a finite mean 
and
$E[|X(1)|\te{\gamma X(1)}]<\infty$. 
\end{As}

In the case of negative drift the  limiting distribution is given as follows:
\begin{Thm}\label{thm}
Let As.~\ref{C} hold.
Then 
$\phi(0)\in(0,\infty)$
and the limit 
$\Psi_\infty(u,v) \doteq\lim_{x\uparrow\infty}\Psi_x(u,v)$ 
exists and is given by 
\begin{equation*}
\Psi_\infty(u,v) = \frac{\gamma}{\phi(0)}\int_u^\infty
\ovl\nu(v+z)\WH V_\gamma(z)\td z,\qquad u,v\geq0,
\end{equation*}
where we denote
$\WH V_\gamma(z) \doteq \int_{[0,z]}\te{\gamma (z-y)} \WH V(\td y)$.
\end{Thm}

\noindent {\bf Remarks.} {\rm (i) If $X$ is spectrally positive 
(i.e.
$\nu((-\infty,0))=0$)
with
$\psi(\theta) = \log E[\te{-\theta X(1)}]$ 
and satisfies As.~\ref{C}, 
we have 
$$
\WH V(y)=\frac{y}{\WH \phi'(0)},\qquad
\phi(0) =\frac{\psi'(0)}{\WH \phi'(0)}, \qquad \WH V_\gamma(z) =
(\gamma\WH \phi'(0))^{-1}(\te{\gamma z} - 1),
$$
where 
$\WH \phi'(0)>0$,
$\gamma$ is
the largest root of $\psi(-\theta)=0$
and
the second equality follows from the Wiener-Hopf factorisation
$-\psi(-\theta)=\phi(\theta)\WH \phi(-\theta)$
and equality
$\WH \phi(0)=0$.
The limit distribution $\Psi_\infty$ is given by 
$$
\Psi_\infty(u,v) = \frac{1}{\psi'(0)}\int_{u}^\infty(\te{\gamma z} - 1) \ovl\nu(v+z)\td z.
$$

\noindent (ii) If $X$ is spectrally positive with $E[X(1)]\in(0,\infty)$,
we have the identities (see e.g.
\cite[p. 191]{Bertoin})
$$\phi'(0) = - \frac{\psi'(0)}{\Phi(0)},
\quad \WH V(y) = \frac{1 - \te{-\Phi(0)y}}{\Phi(0)},
$$
where 
$\psi(\theta)=\log E[\te{-\theta X(1)}]$ 
and 
$\Phi(0)=\WH \phi(0)>0$ 
is the largest root of the equation $\psi(\theta)=0$. 
The joint asymptotic distribution of under- and overshoot is 
in this case given explicitly by the formula
\begin{equation*}
\Psi_\infty(u,v) =
 - \frac{1}{\psi'(0)} \int_u^\infty \ovl\nu(v+y)(1 - \te{-\Phi(0) y})\td y.
\end{equation*}

\noindent (iii) In Corollary~2(ii) of~\cite{MP2},
the marginal law
of
$Z(\infty)$
was identified
and the following expression for the overshoot was given:
\begin{eqnarray}
\label{eq:From_Long_Paper}
P[Z(\infty) > v] = \frac{\gamma}{\phi(0)} \te{-\gamma v}\int_{(v,\infty)} \te{\gamma z}\ovl\nu_H(z)\td z, \q\text{ $v\geq0$},
\end{eqnarray}
where 
$\ovl\nu_H(a) \doteq \nu_H((a,\infty))$, $a>0$,
is the tail of the L\'{e}vy measure $\nu_H$ of $H$.
Combining this with Theorem~\ref{thm}, 
we find that the equality
$\Psi_\infty(0,v)=P[Z(\infty) > v]$
holds for all
$v\geq0$. Indeed,
\begin{eqnarray*}
\frac{\phi(0)}{\gamma} 
\Psi_\infty(0,v) 
& = & 
\int_{[0,\infty)} \WH V(\td y)  \int_y^\infty \te{\gamma (z-y)}\ovl\nu(v+z) \td z = 
\int_{[0,\infty)} \WH V(\td y)  \int_0^\infty \te{\gamma z}\ovl\nu(v+z+y) \td z\\
&= & 
\int_0^\infty \te{\gamma z} \td z
\int_{[0,\infty)} \ovl\nu(v+z+y) \WH V(\td y)  = 
\int_v^\infty \te{\gamma (z-v)} \td z
\int_{[0,\infty)} \ovl\nu(z+y) \WH V(\td y),  
\end{eqnarray*}
which is equal to 
$\frac{\phi(0)}{\gamma}P[Z(\infty) > v]$
by~\eqref{eq:From_Long_Paper} and 
Vigon's identity~\eqref{vi} (established in \cite{Vigon} 
and relating the tail 
$\ovl\nu_H$
of the L\'{e}vy measure $\nu_H$ of $H$ to the dual
renewal function $\WH V$ 
and the upper tail
$\ovl\nu(a)= \nu((a,\infty))$, $a>0$,
of the L\'{e}vy measure~$\nu$ of $X$): 
\begin{equation}\label{vi}
\ovl\nu_H(a) = \int_{[0,\infty)}\ovl\nu(a+y) \WH V(\td y).
\end{equation}
In~\eqref{vi}
the local times
$L$
and
$\WH L$
are normalised so that 
$-\log E[\te{\theta X(1)}]=\phi(-\theta)\WH \phi(\theta)$
holds
for 
$\theta\in[0,\gamma]$
(see e.g.~\cite[Thm.~2.1]{Vigon} and the remark that follows the theorem),
as is assumed throughout this paper.

\noindent (iv) The assumption (in Theorems~\ref{thm:P}
and~\ref{thm})
that 
$X(1)$
is non-lattice is satisfied if the L\'evy measure 
$\nu$
of
$X$
is non-lattice 
or if either the drift or the Gaussian coefficient of
$X$
are non-zero.

\section{Proofs}\label{sec:proofs}
\subsection{Setting} We next describe the setting of
the remainder of the paper, and refer to \cite[Ch. I]{Bertoin} for
further background on L\'{e}vy processes. Let 
$(\Omega,\mc F,\{\mc F(t)\}_{t\ge 0}, P)$ 
be a filtered probability space that carries
a L\'{e}vy process $X$. The sample space $\Omega \doteq D(\mbb R)$
is taken to be the Skorokhod space of real-valued functions that
are right-continuous on $\mbb R_+$ and have left-limits on
$(0,\infty)$, $\{\mc F(t)\}_{t\ge 0}$ denotes the completed
filtration generated by $X$, which is right-continuous, and $\mc
F$ is the completed sigma-algebra generated by 
$\{X(t), t\ge 0\}$. 
For any $x\in\mbb R$ denote by $P_x$ the probability measure
on $(\Omega,\mc F)$ under which the shifted process $X-x$ is a
L\'{e}vy process 
and by 
$E_x$ the expectation under $P_x$. 
Throughout we identify
$P\equiv P_0$
and 
$E\equiv E_0$
and let
$\I_A$ denote the indicator of a set $A$. 

\subsection{Undershoots and overshoots of $X$}
\label{subsec:Under_Over}
An important step in the proof of 
Theorems~\ref{thm:P} and~\ref{thm} consists 
in the identification of the limiting joint distribution of 
the under- and overshoot of $X$, given in Lemma~\ref{lem:Xconv} below.
Let  $T(x) \doteq \inf\{t\ge 0: X(t)> x\}$ denote the first-passage
time of $X$ over the level $x$.  On the set $\{T(x)<\infty\}$,
the overshoot $K(x)$ (resp. undershoot
$k(x)$) of the process 
$X$ at the level $x$ 
is the distance between 
$x$ 
and the positions of $X$ at 
(resp. just before)
the moment $T(x)$: 
\begin{equation}
\label{eq:K_and_k_Levy}
k(x) \doteq x - X(T(x)-),
\qquad K(x) \doteq X(T(x)) - x.
\end{equation}
For any $x>0$ the joint (complementary) distribution of the pair $(k(x),K(x))$ is
denoted by $\Phi_x$, viz.
\begin{eqnarray*}
\Phi_{x}(u,v) \doteq P(k(x) > u, K(x) > v, T(x)<\infty), \qquad u,v\geq0.
\end{eqnarray*}
In the case $P(T(x)<\infty)<1$, or equivalently, when $X$ tends
to $-\infty$ (as is the case under As.~\ref{C}), the  distribution function 
$\Phi_x$  is defective for any $x>0$.
Assuming 
$P(T(x)<\infty)>0$ (as is the case under As.~\ref{C}), 
consider the conditioned distribution
$\Phi_x^\#$
defined as follows:
\begin{equation*}
\Phi_x^\#(u,v) \doteq P(k(x) > u, K(x) > v|T(x)<\infty).
\end{equation*}

\begin{Lemma}\label{lem:Xconv}
(i) Recall that 
$\phi$
is defined in~\eqref{eq:Def_phi}.
Then,
if $X(1)$ is integrable 
with $E[X(1)]\in(0,\infty)$,
it holds
$\phi'(0)=E[H(1)] \in(0, \infty)$.

\noindent (ii)  Suppose $E[|X(1)|]<\infty$ and $E[X(1)]\in(0,\infty)$
and let the law of
$X(1)$
be non-lattice.
Then the
limit $\Phi_\infty(u,v)\doteq\lim_{x\to\infty}\Phi_x(u,v)$ exists
and is given as follows:
\begin{equation}\label{eq:PI}
\Phi_\infty(u,v) = \frac{1}{\phi'(0)}\int_u^\infty \ovl\nu(v+z)\WH
V(z)\td z,\qquad u,v\geq0.
\end{equation}

\noindent (iii) Assume that As. \ref{C} is satisfied. Then
$\Phi^\#_\infty(u,v)
\doteq\lim_{x\to\infty}\Phi^\#_x(u,v)$ exists
and is equal to
\begin{equation}
\label{eq:Cond_Lim}
\Phi^\#_\infty(u,v) = \frac{\gamma}{\phi(0)}\int_u^\infty
\ovl\nu(v+z)\WH V_\gamma(z)\td z,\qquad u,v\geq0,
\end{equation}
where $\WH V_\gamma(z) =\int_{[0,z]}\te{\gamma (z-y)} \WH V(\td y)$
is as defined in Theorem~\ref{thm}.
\end{Lemma}

\noindent {\bf Remark.} The marginal asymptotic distributions of the overshoot
and undershoot of $X$ over a positive level 
under As.~\ref{C} (cf part (iii) of the lemma) were derived 
in~\cite[Thm. 4.1]{GMS}. 
A direct proof of the existence of the joint limit law 
$(k(\infty),K(\infty))$ 
and its explicit description in Lemma~\ref{lem:Xconv}(iii) 
is given in the appendix.


\subsection{Proofs of Theorems~\ref{thm:P} and~\ref{thm}}
In this section we establish our main results.

\subsubsection{Proof of Theorem~\ref{thm:P}:} 
\label{subsec:Proof_of_thm_1}
Fix $0<M<x$ and $u,v\geq0$.
Define
$\WH T(y)\doteq\inf\{t\ge 0: X(t)<-y\}$
for any
$y\geq0$.
The strong Markov property of $Y$ implies the following:
\begin{eqnarray*}
\Psi_x(u,v) & =& 
E[\I_{\{Y(\tau(M))<x\}}P[A_{u,v}(x), \tau(x) < \tau_0|Y({\tau(M)})]]\\
& + &
P[Y(\tau(M))\geq x, A_{u,v}(x)]+
E[P[A_{u,v}(x), \tau(x) > \tau_0|Y({\tau(M)})]],
\end{eqnarray*}
where $A_{u,v}(x) \doteq \{z(x)>u, Z(x) > v\}$ and $\tau_0\doteq\inf\{t\ge0: Y(t)=0\}$. 
Since, for any $y>0$,
the processes $\{Y(t), Y(0)=y, t\leq \tau_0\}$ 
and
$\{X(t), X(0)=y, t\leq \WH T(0)\}$
are equal in law,
for any 
$z\in[M,x)$
we have:
\begin{eqnarray*}
P[A_{u,v}(x), \tau(x) < \tau_0|Y({\tau(M)})=z] 
&=&
\Phi_{x-z}(u,v) - P_z[B_{u,v}(x), \WH T(0) < T(x)],\\
P[A_{u,v}(x), \tau(x) > \tau_0|Y({\tau(M)})=z]
&\leq & P_z[\WH T(0) < T(x)] \leq 
P_z[\WH T(0) < \infty]\leq
P_M[\WH T(0) < \infty],\\
P_z[B_{u,v}(x), \WH T(0) < T(x)] & \leq & P_z[\WH T(0) < T(x)]
 \leq P[\WH T(z) < \infty] \leq P[\WH T(M) < \infty],
\end{eqnarray*}
where we denote
$P_z[\cdot]=P[\cdot|X_0=z]$
and
$B_{u,v}(x) = \{k(x)>u, K(x)>v, T(x)<\infty\}$
(see~\eqref{eq:K_and_k_Levy}).
Since
$\{Y(\tau(M))\geq x, A_{u,v}(x)\}\subset\{\tau(M)=\tau(x)\}$
and
$P[\WH T(M)<\infty]=P_M[\WH T(0)<\infty]$, 
the inequalities above yield the following estimate: 
\begin{eqnarray*}
|\Psi_x(u,v) - \Phi_x(u,v)|  & \leq &   2 P[\WH T(M)<\infty] + 
P[\tau(M)=\tau(x)] \\
& + & \left|E\left[\I_{\{Y(\tau(M))< x\}}\Phi_{x-Y({\tau(M)})}(u,v)\right]- \Phi_x(u,v)\right|.
\end{eqnarray*}
Lemma~\ref{lem:Xconv}~(ii) implies
$\I_{\{Y(\tau(M))< x\}}\Phi_{x-Y({\tau(M)})}(u,v)\to \Phi_\infty(u,v)$
$P$-a.s. 
as
$x\to\infty$.
Also
$\I_{\{\tau(M)=\tau(x)\}}\to 0$
$P$-a.s.
as
$x\to\infty$
and
the dominated convergence theorem yields the following estimate 
for any 
$M>0$:
$$\limsup_{x\to\infty}|\Psi_x(u,v) - \Phi_x(u,v)|\leq 2 P[\WH T(M)<\infty].$$
Since
$E[X(1)]>0$,
$X$ drifts to 
$+\infty$.
Hence
$P[\WH T(M)<\infty]\to 0$
as
$M\to\infty$
and the theorem follows.~\exit 

\subsubsection{Proof of Theorem~\ref{thm}:}
In this section we assume throughout that As. \ref{C} is
satisfied. The proof is based on It\^{o}-excursion theory. Refer
to~\cite[Chs~O,~IV]{Bertoin} for a treatment of It\^{o}-excursion
theory for L\'{e}vy processes and for further references.

Denote by $\e=\{\e(t), t\ge 0\}$ the excursion process of $Y$ away
from zero. Since, under As. \ref{C}, $Y$ is a recurrent strong
Markov process under $P$, It\^{o}'s characterisation yields that
$\e$ is a Poisson point process under $P$. Its intensity measure
under $P$ will be denoted by $n$.
Let $\zeta(\varepsilon)$ denote the lifetime of a generic excursion
$\ve$ and let $\rho(x,\ve)$ denote the first time that the
excursion $\ve$ enters $(x,\infty)$, viz.
\begin{equation}\label{eq:rho}
\rho(x, \varepsilon) = \inf\{t\ge0: \ve(t) > x\}.
\end{equation}
In the sequel we will drop the dependence of $\rho(x,
\varepsilon)$ and $\zeta(\varepsilon)$ on $\varepsilon$, and write
$\zeta$ and $\rho(x)$, respectively.

Theorem \ref{thm} follows directly by combining 
Lemmas~\ref{lem:cond} and~\ref{lem:nasymp} below.

\begin{Lemma}\label{lem:cond}
For any $u,v\geq0$ and $x>0$ the following holds true:
\begin{equation}\label{eq:np}
P(z(x) > u, Z(x) > v) = n(E_{u,v}(x)|\rho(x) < \zeta) \doteq
\frac{n(E_{u,v}(x),\rho(x) < \zeta)}{n(\rho(x) < \zeta)},
\end{equation}
where
$E_{u,v}(x)=\{x-\varepsilon(\rho(x)-)>u,\varepsilon(\rho(x))-x>v,
\rho(x)<\zeta \}$.
\end{Lemma}

\noindent {\it Proof of Lemma \ref{lem:cond}:} 
By~\cite[Ch.~O,~Prop.~2]{Bertoin},
for sets
$A,B$
with
$n(A)\in(0,\infty)$,
we have
$P(\epsilon(T_A)\in B)=n(B|A)=n(A\cap B)/n(A)$
where
$T_A=\inf\{t\geq0:\epsilon(t)\in A\}$.
The lemma now follows by noting that
the left-hand side of~\eqref{eq:np} 
is the probability that the first excursion in
$A=\{\rho(x)<\zeta\}$
is in 
$B=E_{u,v}(x)$.\exit

\begin{Lemma}
\label{lem:nasymp} Let $u,v\geq0$ and recall 
$\WH V_\gamma(z) = \int_{[0,z]}\te{\gamma (z-y)} \WH V(\td y)$.
The following holds true: 
$$
\lim_{x\to\infty}n(E_{u,v}(x)|\rho(x) < \zeta) = \frac{\gamma}{\phi(0)}
\int_u^\infty \ovl\nu(v+z) \WH V_\gamma(z)\td z.
$$
\end{Lemma}

\noindent {\it Remarks.} (i) The proof of Lemma~\ref{lem:nasymp} 
uses the following facts, 
which hold
by~\cite{BertoinDoney}
and~\cite{DoneyMaller},
respectively,
if 
$0$
is regular for 
$(0,\infty)$
under the law of 
$X$
and As.~\ref{C} is satisfied: 
\begin{eqnarray}\label{eq:Cest}
P(T(x) < \infty) & \sim &  C_\gamma\te{-\gamma x} \q\text{as
$x\to\infty$,}\quad\text{ where}\q C_\gamma\doteq 
\frac{\phi(0)}{\gamma\phi'(-\gamma)},\\
n(\rho(x)<\zeta) & \sim & C_\gamma \WH\phi(\gamma) \te{-\gamma x}
\qquad \text{as $x\to\infty$}.
\label{eq:nest}
\end{eqnarray}
Here and throughout the paper we write $f(x)\sim g(x)$ as
$x\to\infty$ if $\lim_{x\to\infty}f(x)/g(x)=1$.


\noindent (ii) A further ingredient of the proof of Lemma~\ref{lem:nasymp} 
are the following asymptotic identities,
established in~\cite[Lemma~9]{MP2}:
\begin{eqnarray}\label{eq:nest2}
n(\te{\gamma \ve(\rho(x))}\I_{\{\rho(x)<\zeta\}}) &\sim&
\WH\phi(\gamma) \qquad \text{as $x\to\infty$},\\
\te{\gamma x}n(\ve(\rho(z)) > x, \rho(x)<\zeta) &=& o(1)\qquad
\text{as $x\to\infty$, for any $z>0$}. \label{eq:nest3}
\end{eqnarray}
(iii) The key observation in~\cite{BertoinDoney}
is that 
$V(\td z)$
is a renewal measure corresponding to the distribution
$P[H(\Theta)\in\td z]$, where
$\Theta$ is an exponential random variable with
$E[\Theta]=1$, 
independent of 
$X$
(and hence of $H$).
Estimate~\eqref{eq:Cest} 
then 
follows from the classical renewal theorem for non-lattice random walks
with the step-size distribution
$H(\Theta)$,
which needs to be non-lattice for the theorem to be applicable
(see the conclusion of the proof of the Theorem in~\cite{BertoinDoney}
for this argument and~\cite[p.~363]{Feller} for the statement of 
the renewal theorem). 
The assumption in~\cite{BertoinDoney},
which ensures this,
stipulates that 
$0$ is regular for
$(0,\infty)$ under the law of
$X$.
Note that this assumption  also 
implies  the non-lattice condition
As.~\ref{C}.
Furthermore,
if 
$X(1)$
is non-lattice,
so is 
$H(\Theta)$
(indeed, if 
$H(\Theta)$
were lattice,
Theorem~30.10 in~\cite{Sato}
would yield that 
$H$
is a compound Poisson process, necessarily with 
a L\'evy measure that has lattice support,
hence implying that 
$X$
itself is a compound Poisson process 
with a L\'evy measure that has lattice support). 
Since the argument
in~\cite{BertoinDoney} 
only requires 
$H(\Theta)$
to be non-lattice,
the estimate 
in~\eqref{eq:Cest} 
remains valid under 
As.~\ref{C}.
Thus the estimate in~\eqref{eq:Cest} holds in our setting. 
Likewise, the argument in~\cite{DoneyMaller}
relies solely on the fact that 
$V(\td z)$
is a renewal measure of a non-lattice law 
and therefore 
estimate~\eqref{eq:nest} also holds under As.~\ref{C}.

\smallskip

\noindent {\it Proof of Lemma \ref{lem:nasymp}:} Fix $M>0$ and pick $u, v\geq0$. The proof
starts from the elementary observation that relates the
following two conditional $n$-measures:
\begin{equation}\label{eq:nc}
n(E_{u,v}(x)|\rho(x)<\zeta) = n(E_{u,v}(x)|\rho(M)<\zeta) \cdot
\frac{n(\rho(M)<\zeta)}{n(\rho(x)<\zeta)},\qquad x>M.
\end{equation}
Recall that the coordinate process under the probability measure
$n(\cdot|\rho(M)<\zeta)$ has the same law as the first excursion
of $Y$ away from zero with height larger than $M$. The strong
Markov property under $n(\cdot|\rho(M)<\zeta)$ implies that
$\varepsilon\circ\theta_{\rho(M)}$ has the same law under
$n(\cdot|\rho(M)<\zeta)$ as the coordinate process of $X$ under $P$,
with entrance law $\mu_M(\td y) \doteq
n(\varepsilon({\rho(M)})\in\td y|\rho(M)<\zeta)$,
that is
killed upon its first entrance into 
$(-\infty,0)$. 
Recall $\WH T(y)=\inf\{t\ge 0: X(t)<-y\}$, for $y\geq0$, 
and note 
that 
for every $x>M$ we have:
\begin{eqnarray}\label{eq:nid}
n(E_{u,v}(x)|\rho(M)<\zeta) =  \int_{[M,x]}
P_z[B_{u,v}(x)] \mu_M(\td z)+O_{u,v}(x),
\end{eqnarray}
where 
$B_{u,v}(x) = \{k(x)>u, K(x)>v, T(x)<\infty\}$
and
$O_{u,v}(x)$ is given by the following expression:
\begin{eqnarray*}
O_{u,v}(x) =  
n(E_{u,v}(x), \varepsilon(\rho(M))>x|\rho(M)<\zeta)-
\int_{[M,x]}
P_z[B_{u,v}(x), \WH T(0)<T(x)<\infty] \mu_M(\td z).
\end{eqnarray*}
Note that
$ P_z[B_{u,v}(x), \WH T(0)<T(x)<\infty]\leq P_z[\WH T(0)<T(x)<\infty]$
for any
$z\in(0,x)$
and hence we find by~\cite[Prop.~7~(i)]{MP2}
that
the following holds (the constant 
$C_\gamma$
is given in~\eqref{eq:Cest}):
\begin{eqnarray*}
P_z[\WH T(0)<T(x)<\infty] & \leq & 
P[\WH T(z)<T(x-z)<\infty]
 \sim  C_\gamma \te{-\gamma x} E[\te{\gamma (X(\WH T(z))+z)}]\qquad\text{as $x\to\infty$.}
\end{eqnarray*}
The following facts hold:
$X(\WH T(z))+z\leq0$,
the measure 
$\mu_M(\td y)$
is concentrated on 
$[M,\infty)$
with 
$\mu_M([M,\infty))=1$
for any $M>0$
and 
equality~\eqref{eq:nest3} 
is satisfied.
Hence, for a fixed
$M>0$,
we have
\begin{equation}\label{eq:n0}
-C_\gamma\te{-\gamma x} \leq O_{u,v}(x)\leq   \te{-\gamma x}o(1)
\qquad \text{as $x\to\infty$.}
\end{equation}

By Lemma~\ref{lem:Xconv}~(iii) we have 
$\lim_{x\to\infty} P_z[k(x)>u, K(x)>v|T(x)<\infty]=\Phi^\#_\infty(u,v)$
for any fixed
$z\geq0$.
Equality~\eqref{eq:Cest} implies
$P_z[T(x)<\infty]=P[T(x-z)<\infty]= \te{-\gamma x}C_\gamma \te{\gamma z}(1+r(x-z))$
as
$x\to\infty$
for any
$z\geq0$,
where 
$r:\RR_+\to\RR$
is bounded and measurable with 
$\lim_{x'\to\infty}r(x')=0$.
By~\eqref{eq:nest2}
$z\mapsto \te{\gamma z}$,
$z\in[M,\infty)$,
is in 
$L^1(\mu_M)$
for all large  
$M$. 
The dominated convergence theorem and~\eqref{eq:nest} 
therefore imply: 
\begin{eqnarray}
\nn
\lim_{x\to\infty}
\int_{[M,x]} \frac{P_z[B_{u,v}(x)]}{n(\rho(x)<\zeta)} \mu_M(\td z)
& = & 
\lim_{x\to\infty}
\int_{[M,x]} P_z[B_{u,v}(x)|T(x)<\infty]\frac{P_z[T(x)<\infty]}{n(\rho(x)<\zeta)} 
\mu_M(\td z)\\
&=&
n(\te{\gamma\ve(\rho(M))}|\rho(M)<\zeta) \cdot
\Phi_\infty^\#(u,v)\cdot \WH \phi (\gamma)^{-1}.
\label{eq:nph}
\end{eqnarray}
Recall that 
$E[X(1)]<0$ 
by As.~\ref{C}
and hence 
$\WH \phi(\gamma)>0$
since 
$\WH H$
is a non-trivial subordinator.
Hence~\eqref{eq:nest}, \eqref{eq:nc}, \eqref{eq:nid}, \eqref{eq:n0} and~\eqref{eq:nph} 
imply the following inequalities 
for any fixed $M>0$:
\begin{eqnarray*}\nn
-\WH \phi (\gamma)^{-1} 
n(\rho(M)<\zeta)
& \leq &
\liminf_{x\to\infty} n(E_{u,v}(x)|\rho(x)<\zeta)- 
\Phi_\infty^\#(u,v) \cdot n(\te{\gamma\ve(\rho(M))} \I_{\{\rho(M)<\zeta\}})
\cdot\WH \phi (\gamma)^{-1}\\
& \leq & 
\limsup_{x\to\infty} n(E_{u,v}(x)|\rho(x)<\zeta)- 
\Phi_\infty^\#(u,v) \cdot n(\te{\gamma\ve(\rho(M))} \I_{\{\rho(M)<\zeta\}})
\cdot\WH \phi (\gamma)^{-1}
\leq 0.
\end{eqnarray*}
Since these inequalities hold for all large $M>0$, 
in the limit as
$M\to\infty$
equation~\eqref{eq:nest2} implies
$\lim_{x\to\infty}n(E_{u,v}(x)|\rho(x)<\zeta)
=\Phi_\infty^\#(u,v)$.
This, together with
Lemma~\ref{lem:Xconv}~(iii), 
concludes the proof.\exit

\appendix
\section{Proof of Lemma~\ref{lem:Xconv}}
\label{app:proof_lemma}
\begin{proof}
(i) $E[X(1)]>0$ implies 
$X(t)\to\infty$ $P$-a.s. as $t\uparrow\infty$.
Hence
$E[H(1)]\in(0,\infty]$
and
$E[\WH L(\infty)]=1/\WH \phi(0)<\infty$.
Since $E[X(1)]<\infty$, we have
$\int_{[1,\infty)}y\nu(\td y)<\infty$.
By definition 
we have
$\WH V(\infty)\doteq\lim_{y\to\infty}\WH V(y)=E[\WH L(\infty)]$.
Fubini's theorem,
the estimate 
$\int_{[1,\infty)}z\nu(y+\td z)\leq
\int_{[1,\infty)}(z+y)\nu(y+\td z)\leq
\int_{[1,\infty)}x\nu(\td x)<\infty$
for any
$y\geq0$
and~\eqref{vi} imply
\begin{eqnarray*}
\int_{[1,\infty)}y\nu_H(\td y) & = & 
\int_{[0,\infty)} \WH V(\td y)
\int_{[1,\infty)}z\nu(y+\td z)\leq
\WH V(\infty)
\int_{[1,\infty)}x\nu(\td x)<\infty, 
\end{eqnarray*}
and part~(i) of the lemma follows.
 
(ii) The compensation formula applied to the Poisson
point process $\{\Delta X(t),t\geq0\}$ 
(here $X(0-)\doteq0$
and
$\Delta X(t)\doteq X(t)-X(t-)$ for 
$t\geq0$)
and the form of the resolvent of $X$
killed upon entering $(x,\infty)$ (see \cite[p.176]{Bertoin}) 
imply the following 
identity
(recall 
$\ovl\nu(a)=\nu((a,\infty))$, 
$a>0$,
is the tail of the L\'evy measure 
$\nu$):
\begin{eqnarray}
\nonumber
\Phi_x(u,v) & = & E\sum_{t>0}\I_{\{X^*(t-)<x,\>\> x-X(t-)>u,\>\>X(t-)+\Delta X(t)-x>v\}}\\
\nonumber
& = & E\int_0^\infty\ovl \nu(v+x-X(t-))\I_{\{x-X(t-)>u,X^*(t-)<x\}}\td t
=E\int_0^{T(x)}\ovl \nu(v+x-X(t))\I_{\{x-X(t)>u\}}\td t\\
& = & \int_{[0,x]}\ovl F(x-z)  V(\td z), 
\label{eq:Phix}
\end{eqnarray}
where $\ovl F$ is given by the expression 
\begin{equation}
\label{eq:ovl_F}
\ovl F(z) =
\int_{[0,\infty)}\ovl\nu(v+z+y)\I_{\{z+y>u\}}
\WH{V}(\td y),\qquad\text{for any $z\geq0$,}
\end{equation}
and the function 
$V$
is defined in~\eqref{eq:LTV}
(an argument based on the quintuple law 
from~\cite[Thm.~3]{DK}
can also be applied to establish~\eqref{eq:Phix}).
The function
$z\mapsto \ovl\nu(z+v)$
is integrable on 
$(0,\infty)$
by the assumption that
$E[|X(1)|]<\infty$
and the inequality 
$\WH V(\infty)<\infty$
holds 
(see e.g. proof of Lemma~\ref{lem:Xconv}(i) above).
Hence the inequalities 
$0\leq \ovl F(z)\leq \ovl\nu(z+v)\WH V(\infty)$
hold for all
$z\geq0$,
making the function 
$\ovl F$
directly Riemann integrable 
as defined in~\cite[Definition on p.~362]{Feller}.

Let $\Theta$ be  
independent of 
$H$
and exponentially distributed with
$E[\Theta]=1$. 
The law
$P[H(\Theta)\in\td z]$ 
has the mean equal to 
$\phi'(0)$.
By Remark~(iii) following Lemma~\ref{lem:nasymp},
the renewal theorem in~\cite[Thm.~on~p.~363]{Feller}
and~\eqref{eq:Phix}
imply that 
$\Phi_\infty(u,v)=\lim_{x\to\infty}\Phi_x(u,v)$
exists and is equal to 
$$
\Phi_\infty(u,v)=\frac{1}{\phi'(0)}\int_{[0,\infty)}\ovl F(z)\td z. 
$$

\noindent
The definition of 
$\ovl F$
in~\eqref{eq:ovl_F}
and several applications of Fubini's theorem 
yield the following:
\begin{eqnarray}
\nonumber
\int_{[0,\infty)}\ovl F(z)\td z &=&
\int_{[0,u]}\td z \int_{(u-z,\infty)} \ovl \nu(v+z+y)\WH V(\td y) 
+
\int_{(u,\infty)}\td z \int_{[0,\infty)} \ovl \nu(v+z+y)\WH V(\td y) \\
\nonumber
& = &
\int_{[0,\infty)} \WH V(\td y)
 \int_{((u-y)^+,u]} \ovl \nu(v+z+y)\td z
+
\int_{[0,\infty)} \WH V(\td y)
 \int_{(u,\infty)} \ovl \nu(v+z+y)\td z \\
& = & 
\int_{[0,\infty)}\WH V(\td y)\int_{[(u-y)^+,\infty)}\ovl \nu(v+z+y)\td z,
\label{eq:Final_Eq_F_ovl}
\end{eqnarray}
where as usual
$(u-y)^+=\max\{u-y,0\}$.
The equality in~\eqref{eq:Final_Eq_F_ovl}
and further applications of the Fubini theorem
imply part~(ii) of the lemma:
\begin{eqnarray*}
\int_{[0,\infty)}\ovl F(z)\td z &=&
\int_{[0,u]}\WH V(\td y)\int_{[u-y,\infty)}\ovl \nu(v+z+y)\td z
+
\int_{(u,\infty)}\WH V(\td y)\int_{[0,\infty)}\ovl \nu(v+z+y)\td z
\\
& = & 
\int_{[0,u]}\WH V(\td y)\int_{[u,\infty)}\ovl \nu(v+z)\td z
+
\int_{(u,\infty)}\WH V(\td y)\int_{[y,\infty)}\ovl \nu(v+z)\td z\\
& = &
\WH V(u) \int_{[u,\infty)}\ovl \nu(v+z)\td z + 
\int_{[u,\infty)}\ovl \nu(v+z)\td z
\int_{(u,z]} \WH V(\td y)\\
& = & 
\WH V(u) \int_{[u,\infty)}\ovl \nu(v+z)\td z + 
\int_{[u,\infty)}\left(\WH V(z)- \WH V(u)\right) \ovl \nu(v+z)\td z\\
& = & \int_{[u,\infty)}\WH V(z)\ovl \nu(v+z)\td z.
\end{eqnarray*}

(iii) 
Let $P^{(\gamma)}$ be the Cram\'{e}r measure
on $(\Omega,\mc F)$, the restriction of
which to $\mc F(t)$ 
is defined by
$P^{(\gamma)}(A) \doteq E[\te{\gamma X(t)}\I_A]$
for any 
$A\in\mc
F(t),t\in\mbb R_+.$
Under
$P^{(\gamma)}$ 
it holds
$E^{(\gamma)}[|X(1)|]=E[|X(1)|\te{\gamma X(1)}]<\infty$
and
$E^{(\gamma)}[X(1)]>0$ and hence
$P^{(\gamma)}(T(x)<\infty)=1$.
Define
$\Phi_{x}^{(\gamma)}(u,v) \doteq P^{(\gamma)}(k(x) > u, K(x) > v, T(x)<\infty)$
for any 
$u,v\geq0$.
Changing the measure yields
$$
\Phi^\#_x(u,v)P(T(x)<\infty) = \te{-\gamma x}
E^{(\gamma)}[\te{-\gamma K(x)}\I_{\{k(x)>u, K(x)>v, T(x)<\infty\}}]=
\te{-\gamma x} \int_{(v,\infty)}\te{-\gamma w}\Phi_x^{(\gamma)}(u,\td w).
$$

By part~(ii) of the lemma, the limit 
$\Phi_{x}^{(\gamma)}(u,v) \to\Phi_\infty^{(\gamma)}(u,v)$,
as
$x\uparrow\infty$,
exists 
for all
$u,v\geq0$.
Assume first
$\Phi_\infty^{(\gamma)}(u,v)>0$
and note that 
the probability measures
$\I_{\{w>v\}}\Phi_x^{(\gamma)}(u,\td w)/\Phi_x^{(\gamma)}(u,v)$
on $\RR$ 
converge weakly to the probability measure
$\I_{\{w>v\}}\Phi_\infty^{(\gamma)}(u,\td w)/\Phi_\infty^{(\gamma)}(u,v)$
as
$x\uparrow\infty$.
Hence~\cite[Thm.3.9.1(vi)]{Durrett}
applied to the bounded function
$w\mapsto \I_{\{w>v\}} \te{-w\gamma}$,
the Cram\'{e}r's asymptotics in~\eqref{eq:Cest}
and Lemma~\ref{lem:Xconv}~(ii)
imply the
following equalities 
\begin{equation}
\label{eq:Hopefully_Final}
\lim_{x\to\infty}
\Phi^\#_x(u,v) = C_\g^{-1} \int_{(v,\infty)} \te{-\g
w}\Phi_\infty^{(\gamma)}(u,\td w) = 
\frac{C_\g^{-1}}{\phi^{(\gamma)\prime}(0)}\int_{(v,\infty)}
\te{-\gamma w}\int_{(u,\infty)}\nu^{(\gamma)}(y+\td w) \WH V^{(\gamma)}(y)\td y.
\end{equation}
In the case
$\Phi_\infty^{(\gamma)}(u,v)=0$
we note
$\int_{(v,\infty)}\te{-\gamma w}\Phi_x^{(\gamma)}(u,\td w)\leq 
\Phi_x^{(\gamma)}(u,v)$.
Hence by~\eqref{eq:Cest} and Lemma~\ref{lem:Xconv}~(ii)
we find
$ \lim_{x\to\infty} \Phi^\#_x(u,v) = 
\lim_{x\to\infty}\frac{\te{-\gamma x}}{P(T(x)<\infty)}\int_{(v,\infty)}\te{-\gamma w}\Phi_x^{(\gamma)}(u,\td w)=0$.
Therefore the first equality in~\eqref{eq:Hopefully_Final} holds also in the case
$\Phi_\infty^{(\gamma)}(u,v)=0$.

The Wiener-Hopf factorisation~\cite[p.~166, Eqn.~(4)]{Bertoin}
implies
$\phi^{(\gamma)}(\theta)=\phi(\theta-\gamma)$ 
and
$\WH \phi^{(\gamma)}(\theta)=\WH \phi(\theta+\gamma)$ 
for all
$\theta\geq0$.
The elementary equality
$\nu^{(\gamma)}(\td y) = \te{\gamma y}\nu(\td y)$ 
and the form of $C_\gamma$ given in~\eqref{eq:Cest} 
therefore yield 
$$
\lim_{x\to\infty} 
\Phi^\#_x(u,v) = \frac{\gamma}{\phi(0)}\int_u^\infty
\ovl\nu(v+y)\te{\gamma y} \WH V^{(\gamma)}(y)\td y. 
$$

The Laplace transform of 
$\WH V^{(\gamma)}$
is given by
$[\theta\WH\phi^{(\gamma)}(\theta)]^{-1}=[\theta\WH\phi(\theta+\gamma)]^{-1}$ 
(cf.~\eqref{eq:LTV}).
It follows that 
the Laplace transforms
of the function
$y\mapsto\te{\gamma y} \WH V^{(\gamma)}(y)$ 
and the convolution
$y\mapsto \WH V_\gamma(y)=\int_{[0,y]}\te{\gamma(y-z)} \WH V(\td z)$
are both equal to 
$[(\theta-\gamma)\WH\phi(\theta)]^{-1}$
(recall that
$\int_{[0,\infty)}\te{-\theta z}\WH V(\td z)=1/\WH \phi(\theta)$). 
Hence the two functions can only differ on a set with at most countably many
points, which has Lebesgue measure zero.
Therefore the formula in~\eqref{eq:Cond_Lim}
and the lemma follow.
\end{proof}

\end{document}